# Alternative Orthogonal Polynomials

Vladimir Chelyshkov

*Institute of Hydromechanics of the NAS, Ukraine,
Georgia Southern University, USA*

**Abstract.** The double-direction orthogonalization algorithm is applied to construct sequences of polynomials, which are orthogonal over the interval $[0,1]$ with the weighting function 1. Functional and recurrent relations are derived for the sequences that are the result of the inverse orthogonalization procedure. Quadratures, generated by the sequences, are described.

Family of common orthogonal polynomials originates from a problem on the differential equation of the hypergeometric type, which solution is subjected to certain additional requirements [1]. The polynomials also may be defined by an orthogonalization procedure, if it is applied to the fundamental sequence $\{x^k\}$ in the order of power increase. Generalizing this approach, one can develop the orthogonalization procedure beginning with an arbitrary number of the sequence, both in the direct and inverse order. The double-direction algorithm of orthogonalization was introduced in [2] for defining orthogonal sequences of exponents. It was also mentioned there that the algorithm may be applied to the fundamental sequence under various orthogonality relations, and, for the polynomials constructed, the inverse algorithm retains the properties of the original sequence as $x \to 0$, if $x \in [0,1]$. Here we describe an example of such alternative sequences. The alternative orthogonal polynomials (AOP) obtained are not solutions of the equation of the hypergeometric type, but they can be expressed in terms of the Jacobi polynomials.

Let $n$ be a fixed whole number, $\mathcal{P}_n$ and $P_n$ are sequences of polynomials,

$$\mathcal{P}_n = \{\mathcal{P}_{nk}\}_{k=n}^{0}, \quad \mathcal{P}_{nk} \equiv \mathcal{P}_{nk}(x) = \sum_{l=k}^{n} \tau_{nkl} x^l, \tag{1}$$

$$P_n = \{P_{nk}\}_{\infty}^{k=n}, \quad P_{nk} \equiv P_{nk}(x) = \sum_{l=n}^{k} T_{nkl} x^l, \tag{2}$$

that hold orthogonality relationships

$$\int_0^1 \mathcal{P}_{nk} \mathcal{P}_{nl} dx = \begin{cases} 0 & k \neq l \\ 1/(k+l+1) & k = l \end{cases} \quad k,l = 0,1,\ldots,n, \tag{3}$$

$$\int_0^1 P_{nk} P_{nl} dx = \begin{cases} 0 & k \neq l \\ 1/(k+l+1) & k = l \end{cases} \quad k,l = n, n+1,..., \tag{4}$$

and normalizing conditions

$$\operatorname{sign}(\tau_{nkn}) = (-1)^{n-k}, \quad \operatorname{sign}(T_{nkk}) = (-1)^{k-n}. \tag{5}$$

The coefficients $\tau_{nkl}$ and $T_{nkl}$ of the polynomials $\mathcal{P}_{nk}$ and $P_{nl}$ are defined uniquely by requirements (3) - (5) and the Gram-Schmidt orthogonalization algorithm, which is realized in the order of decreasing $k$ from $n$ to 0 for sequences (1), and in the order of increasing $k$, originating from $k = n$, for sequences (2). The sequences $\mathcal{P}_n$ and $P_n$ have different properties, if $x \sim 0$. For fixed $n$ and $x \to 0$ $\mathcal{P}_{nk}(x) \sim x^k$, $k = 0, 1, ... n$, and $P_{nk}(x) \sim x^n$, $k = n, n+1, ...$. The sequence $P_0$ represents the shifted to the interval $[0,1]$ Legendre polynomials; $P_n$, $n > 0$, is considered as an auxiliary sequence, and the sequence $\mathcal{P}_n$ is introduced here as the *alternative Legendre polynomials* (ALP). The polynomials $\mathcal{P}_{nk}$ and $P_{nk}$ have properties, which are analogues to the properties of common orthogonal polynomials. Since $n$ is fixed, the polynomials $P_{nl}(x)$ can be immediately connected to a fixed set of the Jacobi polynomials $P_m^{(\alpha,\beta)}(\xi)$ [1] by verifying property (4) and (5), and the following representation holds

$$P_{nk}(x) = x^n P_{k-n}^{(2n,0)}(1-2x). \tag{6}$$

This relation can be used directly to describe the properties of $P_{nk}$, and one of the formulas that follow from (6) is the integral representation

$$P_{nk}(x) = \frac{1}{2\pi i} \frac{1}{x^n} \int_C \frac{z^{k+n}(1-z)^{k-n}}{(z-x)^{k-n+1}} dz. \tag{7}$$

Here $C$ is a closed curve, which encloses point $z = x$.

Realizing the orthogonalization procedure one can suppose that the explicit definition of the polynomials $\mathcal{P}_{nk}$ is

$$\mathcal{P}_{nk}(x) = \sum_{j=0}^{n-k} (-1)^j \binom{n-k}{j} \binom{n+k+1+j}{n-k} x^{k+j}, \quad k = 0, 1, ..., n. \tag{8}$$

This yields Rodrigues' type representation,

$$\mathcal{P}_{nk}(x) = \frac{1}{(n-k)!} \frac{1}{x^{k+1}} \frac{d^{n-k}}{dx^{n-k}} (x^{n+k+1}(1-x)^{n-k}), \quad k = 0, 1, ..., n, \tag{9}$$

and orthogonality relationship (3) is confirmed by applying last formula. It also follows from (9) that

$$\int_0^1 \mathcal{P}_{nk}(x)dx = \int_0^1 x^n dx = \frac{1}{n+1}. \tag{10}$$

Making use of formula (9), and the integral formula Cauchy for derivatives of an analytic function one can obtain the integral representation

$$\mathcal{P}_{nk}(x) = \frac{1}{2\pi i} \frac{1}{x^{n+2}} \int_{C_1} \frac{z^{-(n+k+2)}(1-z)^{n-k}}{(z-x^{-1})^{n-k+1}} dz, \tag{11}$$

where $C_1$ is a closed curve, which encloses point $z = x^{-1}$. Representations (7) and (11) lead directly to *the reciprocity relation*

$$\mathcal{P}_{nk}(x) = x^{-1} P_{-(n+1),-(k+1)}(x^{-1}). \tag{12}$$

Relationship, similar to (12), holds also for orthogonal exponential polynomials [2].

Formula (12) facilitates description of the ALP, and the results that are shown below can be obtained making use of the auxiliary sequences $P_{nk}(x)$ and relationship (12). In particular,

$$\mathcal{P}_{nn} = x^n, \quad \mathcal{P}_{n,n-1} = 2nx^{n-1} - (2n+1)x^n,$$

and the following recurrence relations and differentiation formulas hold:

$$a_{nk}\mathcal{P}_{n,k-1} = (b_{nk}x^{-1} - c_{nk})\mathcal{P}_{nk} - d_{nk}\mathcal{P}_{n,k+1} \tag{13}$$

$$\alpha_{nk}x(1-x)\mathcal{P}'_{nk} = (\beta_{nk} - \gamma_{nk}x)\mathcal{P}_{nk} - \delta_{nk}x\mathcal{P}_{n,k+1}, \tag{14}$$

$$\kappa_{nk}x(x-1)\mathcal{P}'_{nk} = (\lambda_{nk} - \mu_{nk}x)\mathcal{P}_{nk} - \nu_{nk}x\mathcal{P}_{n,k-1}, \tag{15}$$

where

$$a_{nk} = (k+1)(n-k+1)(n+k+1), \quad b_{nk} = k(2k+1)(2k+2),$$

$$c_{nk} = (2k+1)((n+1)^2 + k^2 + k), \quad d_{nk} = k(n-k)(n+k+2),$$

$$\alpha_{nk} = 2(k+1), \quad \beta_{nk} = 2k(k+1), \quad \gamma_{nk} = n^2 + k^2 + 2n,$$

$$\delta_{nk} = (n-k)(n+k+2),$$

$$\kappa_{nk} = 2k, \quad \lambda_{nk} = 2k(k+1), \quad \mu_{nk} = (n+1)^2 + k^2,$$

$$\nu_{nk} = (n-k+1)(n+k+1). \tag{16}$$

The polynomial $x\mathcal{P}_{nk}(x)$ is a solution of the differential equation

$$x^2(1-x)\zeta'' - x^2\zeta' + ((n+1)^2 x - k(k+1))\zeta = 0, \quad k = 0,1,...,n \tag{17}$$

that also follows from the constructions developed. Making use of the substitution $\varsigma(x) = x^k u(x)$ one can represent the polynomial solution of equation (17) in terms of the hypergeometric function $F$, and the following relationships hold

$$\mathcal{P}_{nk}(x) = \binom{n+k}{n-k} x^k F(k-n, k+n+2; 2k+1; x),$$

$$\mathcal{P}_{nk}(x) = x^k P_{n-k}^{(2k,1)}(1-2x), \quad k = 0,1,...,n.$$

Thus, the ALP are related to *different* families of the Jacoby polynomials.

Let $f(x) \in C[0,1]$. The ALP $\mathcal{P}_{nk}(x) \sim x^k$ as $x \to 0$, and they can be applied for constructing the Gauss-type quadrature formula. Making use of the orthogonality relations and the approach developed in [3, page 378], one can obtain the formulas

$$\int_0^1 f(x)dx \approx \sum_{j=k}^n w_j f(x_j), \tag{18}$$

where $x_j$ are the zeros of the polynomial $\mathcal{P}_{n,k-1}(x)$, $k=1,2,...n$, and the weighting factors $w_j$ are

$$w_j = \frac{1}{\sum_{l=k}^n (2l+1)\mathcal{P}_{nl}^2(x_j)}. \tag{19}$$

Alternative Gauss quadrature (AGQ) (18), (19) is exact for $x^l$, $2k-1 \leq l \leq 2n$. Analogues approach was described for orthogonal exponential polynomials in [4]. Taking into account property (10) of the ALP one can come to the conclusion that the Radau-type quadrature is necessary to add the function $x^{2k-2}$ to the above ones exactly integrable. In particular, if $k=1$, such a quadrature is correct for the function $f(x) \equiv 1$. To avoid application of the alternative Radau quadrature, *almost* orthogonal sequence of polynomials $\{1, \mathcal{P}_{nk}\}_{k=n}^{1}$ coupled with the AGQ might be considered for problems on approximation.

Any polynomial function, say, $P_m(x)$, $x \in [0,1]$ can be considered as a term of the alternative polynomial hierarchy $\mathbb{A} = \{\mathbf{Q}_n\}_{n=0}^{\infty}$, $\mathbf{Q}_n = \{Q_{nk}(x)\}_{k=n}^{0}$, $Q_{nk}(x) = x^k P_{n-k}(x)$ where $n-k$ is the span of the polynomial $Q_{nk}(x)$. Being subjected to certain orthogonality relations, alternative polynomials keep distinctively all the attributes of regular orthogonal polynomials. Thus, peculiarity of such a construction has led to description of the sequence of orthogonal exponential polynomials that generate Gaussian quadrature for exponents [2], [4]. The algorithm of inverse orthogonalization of the fundamental sequence results in redistribution of the zeros of common orthogonal polynomials that makes possible developing different Gauss-type quadratures. This facilitates application of the AOP to integrating the initial value problem [5].